\documentclass[12pt]{article} 
\usepackage{amssymb}
\usepackage{amsthm}
\usepackage{times} 
\usepackage[all]{xy}
\mathversion{bold}
\newcommand{\nc}{\newcommand} 
\nc{\bb}{\bigskip}
\nc{\C}{\mathbb{C}}
\nc{\cl}{\centerline} 
\nc{\dsp}{\displaystyle} 
\nc{\D}{\mathop{\rm D}\nolimits} 
\nc{\End}{\mathop{\rm End}\nolimits} 
\nc{\epi}{\twoheadrightarrow}
\nc{\Ext}{\mathop{\rm Ext}\nolimits} 
\nc{\f}{\varphi}
\nc{\go}{\mathfrak} 
\nc{\Hom}{\mathop{\rm Hom}\nolimits}
\nc{\Id}{\mathop{\rm Id}\nolimits} 
\nc{\ind}{\hskip 1em\relax}
\nc{\K}{\mathbb{K}}
\nc{\LL}{\mathop{\rm L}\nolimits}
\nc{\mod}{\mathop{\rm mod}\nolimits}  
\nc{\mono}{\rightarrowtail} 
\nc{\pt}{\bullet}
\nc{\RR}{\mathop{\rm R}\nolimits} 
\nc{\rad}{\mathop{\rm rad}\nolimits} 
\nc{\then}{\Longrightarrow}
\nc{\Tor}{\mathop{\rm Tor}\nolimits} 
\nc{\Z}{\mathbb{Z}}

\begin{document}
\theoremstyle{definition}
\newtheorem{ttt}{Theorem}
\newtheorem{ddd}[ttt]{Definition}
\newtheorem{lem}[ttt]{Lemma}
\newtheorem{ppp}[ttt]{Proposition}

\cl{\Large\bf About a Theorem of Cline, Parshall and Scott}\bb\bb

Let $\go g$ be a complex semisimple Lie algebra, let 
$\go{h\subset b}$ be respectively Cartan and Borel 
subalgebras of $\go g$, put $\go{n:=[b,b]}$, say that the roots of $\go h$ 
in $\go n$ are positive, let ${\cal W}$ be the Weyl group 
equipped with the Bruhat ordering, let ${\cal O}_0$ be the category of 
those BGG-modules which have the generalized infinitesimal character of the
trivial module. The simple 
objects of ${\cal O}_0$ are parametrized by ${\cal W}$. 
Say that $Y\subset {\cal W}$ is an {\bf initial segment} if $x\le y$ and 
$y\in Y$ imply $x\in Y$, and that $w\in {\cal W}$ lies in the {\bf support} of 
$V\in {\cal O}_0$ if the simple object attached to $w$ is a 
subquotient of $V$. For such an initial segment $Y$ let 
${\cal O}(Y)$ be the subcategory of 
${\cal O}_0$ consisting of objects supported on $Y\subset {\cal W}$, and
let $i:{\cal O}(Y)\to{\cal O}_0$ be the inclusion. Theorem 3.9 of Cline, 
Parshall and Scott in \cite{cps2} implies the following. 

\begin{ttt}\label{cpst} The functor $i_*:\D^b({\cal O}(Y))\to 
 \D^b({\cal O}_0)$ admits a left adjoint $i^*$ and a right adjoint 
 $i^!$ satisfying 
 $i^*\ i_*\simeq \Id_{\D^b({\cal O}(Y))}\simeq i^!\ i_*$. 
 In particular $i_*$ is a full embedding. 
\end{ttt}

\ind The purpose of this text is to give a simple proof of this 
Theorem and to suggest an analog for Harish-Chandra modules. \bb

\ind Theorem \ref{cpst} above will follow from Theorem \ref{Db} below. 
Say that a {\bf BGS category} is an abelian category satisfying 
Conditions (1) to (6) in Section 3.2 of Beilinson, Ginzburg and 
Soergel \cite{bgs}. By Theorem 3.2.1 and Corollary 3.2.2 in \cite{bgs}, 
Theorems \ref{cpst} and \ref{Db} apply to BGS categories. In \cite{bgs} many 
natural examples of BGS categories are given, like 
(in the notation of \cite{bgs}) the categories of BGG modules 
$\cal{O}_\lambda$ and $\cal{O}^\go{q}$ defined in Section 1.1, or more 
generally the category ${\cal P}(X,\cal W)$ of perverse sheaves considered in 
Section 3.3. The fact that ${\cal P}(X,\cal W)$ is a BGS category is 
viewed as obvious in \cite{bgs} (and I have no doubt that it is so 
for algebraic geometers). Theorem 3.5.3 of \cite{bgs} implies that 
$\cal{O}^\go{q}$ is of the form ${\cal P}(X,\cal W)$, and Theorem 3.11.1 
of \cite{bgs} entails that $\cal{O}_\lambda$ is opposite to 
$\cal{O}^\go{q}$ (for some $\go q$). Since the axioms of BGS categories 
are selfopposite, $\cal{O}_\lambda$ is BGS. \bb

\ind Thank you to Bernhard Keller and Wolfgang Soergel for their interest,  
and to Martin Olbrich for having pointed out some mistakes in a previous 
version. 


\section{Statement}

\ind Let $A$ be a ring, $X$ a finite set and $e_\pt=(e_x)_{x\in X}$ a 
family of idempotents of $A$ satisfying 
$\sum_{x\in X}e_x=1$ and $e_xe_y=\delta_{xy}e_x$ (Kronecker delta) 
for all $x,y\in X$. \bb 

\ind The {\bf support} of an $A$-module $V$ is the set 
$\{x\in X\ |\ e_xV\not=0\}$. Let $\le$ be a partial ordering on $X$, and 
for any initial segment $Y$ put 

$$A(Y):=A\left/\sum_{x\notin Y}Ae_xA,\right.$$

so that $A(Y)$-{\bf mod} is the full subcategory of $A$-{\bf mod} whose 
objects are supported on $Y$. (Here and in the sequel, for any ring $B$, 
we denote by $B$-{\bf mod} the category of $B$-modules.) The image 
of $e_y$ in $A(Y)$ will be still denoted by $e_y$. \bb

\ind Assume that, for any pair $(Y,y)$ where $Y$ is an initial segment 
and $y$ a maximal element of $Y$, the module $M_y:=A(Y)e_y$ does {\bf not} 
depend on $Y$, but only on $y$. This is equivalent to the 
requirement that $A(Y)e_y$ be supported on $$\{x\in X\ |\ x\le y\}.$$

\ind If $(V_\gamma)_{\gamma\in\Gamma}$ a family of $A$-modules, let 
$\langle V_\gamma\rangle_{\gamma\in\Gamma}$ 
denote the class of 
those $A$-modules which admit a finite filtration whose 
associated graded object is isomorphic to a product of members of the 
family. \bb

\ind Assume that, for any $x\in X$, the module $Ae_x$ belongs to 
$\langle M_y\rangle_{y\in X}$.  


\begin{ttt} \label{Db} Let $Y\subset X$ be an initial segment and 
 $i_*:\D^b(A(Y)\mbox{-{\bf mod}})\to\D^b(A\mbox{-{\bf mod}})$ the 
 induced functor. Then $i^!:=\RR\Hom_A(A(Y),?)$ is a right adjoint to the 
 functor $i_*$ from $\D^b(A(Y)\mbox{-{\bf mod}})$ to 
 $\D^b(A\mbox{-{\bf mod}})$ and we have 
 $i^!\ i_*\simeq \Id_{\D^b(A(Y)-\mod)}$. In particular $i_*$ is a full 
 embedding. If the right flat dimension of $A(Y)$ over $A$ is finite, then 
 $i^*:=A(Y)\otimes_A^{\LL}?$ is a left adjoint to $i_*$ satisfying 
 $i^*\ i_*\simeq \Id_{\D^b(A(Y)-\mod)}$. 
\end{ttt}


\section{Proof}

{\bf Proof that Theorem \ref{Db} implies Theorem \ref{cpst}.} 
In view of BGG \cite{bgg1} it suffices to check that the Verma module 
$M_x$, with $x\in\cal W$, is projective into ${\cal O}({\cal W}_{\not>x})$ 
(obvious notation). Let $V$ be in ${\cal O}({\cal W}_{\not>x})$ and, 
for any $\lambda\in\go{h}^*$, let $V^\lambda$ be the corresponding weight 
subspace of $V$. Letting $\rho$ be the half sum of the positive roots and 
putting $\lambda:=-x\rho-\rho$ we have 

$$\Hom_\go{g}(M_x,V)\simeq H^0(\go{n},V^\lambda)\subset V^\lambda.$$ 

It suffices to show that this inclusion is an equality. Otherwise 
there would be a weight $\mu$ satisfying

$$\mu>\lambda,\quad V^\mu\not=0,\quad \go{n}V^\mu=0.$$

Letting $L_y$ be a simple quotient of $U(\go{g})V^\mu$, we would have 

$$-y\rho-\rho=\mu>\lambda=-x\rho-\rho$$ 

and thus (see for instance Lemma 7.7.2 in Dixmier \cite{dix}) 
$y>x$, which is impossible. $\square$

\begin{lem}\label{extr} Let $A$ be a ring, $I$ a left projective 
 idempotent twosided ideal and $B$ the quotient ring $A/I$. 
 Then $i^!:=\RR\Hom_A(B,?)$ is a right adjoint to the functor $i_*$ from 
 $\D^b(B\mbox{-{\bf mod}})$ to $\D^b(A\mbox{-{\bf mod}})$ and we have 
 $i^!\ i_*\simeq \Id_{\D^b(B-\mod)}.$ 
 In particular $i_*$ is a full embedding. If the right flat dimension of 
 $B$ over $A$ is finite, then $i^*:=B\otimes_A^{\LL}?$ is a left adjoint to 
 $i_*$ satisfying $i^*\ i_*\simeq \Id_{\D^b(B-\mod)}$. 
\end{lem}

{\bf Proof.} The Lemma follows from Theorem 3.1 and Proposition 3.6 of Cline, 
Parshall and Scott in \cite{cps1}. $\square$\bb

\ind Let us go back to the setting of Theorem~\ref{Db}. 

\begin{lem}\label{swap} For any $x,y\in X$ with $x$ maximal  
there is a nonnegative integer $n$ and an exact sequence 
$(Ae_x)^n\mono Ae_y\epi V$ such that $V\in\langle M_z\rangle_{z<x}$. 
In particular $e_xV=0$. 
\end{lem}

{\bf Proof.} This follows from the projectivity of $M_x=Ae_x$. 
$\square$\bb

{\bf Proof of Theorem~\ref{Db}.} 
Assume $Y=X\backslash \{x\}$ where $x$ is maximal. 
Put $e:=e_x$, $I:=AeA$ and $B:=A(Y)=A/I$. By the previous Lemma 
there is a nonnegative integer $n$ and an exact sequence 
$(Ae)^n\mono A\epi V$ with $IV=0$. Letting $J\subset A$ be the image of 
$(Ae)^n\mono A$, we have $J=IJ\subset I\subset J$, and thus $I=J$. 
Lemma~\ref{extr} applies, proving the Theorem for 
the particular initial segment $Y$. Lemma~\ref{swap} shows that 
$(B,Y,(e_y)_{y\in Y})$ satisfies the assumptions of Theorem~\ref{Db}~; 
and an obvious induction shows the existence of some right adjoint 
$i^!$ to $i_*$, which is a full embedding satisfying 
$i^!\ i_*\simeq\Id_{\D^b(A(Y)-\mod)}$. Then the Theorem follows from
Theorem 3.1 of \cite{cps1}. $\square$


\section{Harish-Chandra modules}

\ind Let $G$ be a connected semisimple Lie group with finite center, 
let $K$ be a maximal compact subgroup, 
let $Z$ be the center of the complexified enveloping algebra, 
let $I$ be the annihilator of the trivial module $\C$ in $Z$, 
let $\hat{Z}$ and $\hat{I}$ be the respective $I$-adic completions of 
$Z$ and $I$, and let ${\cal H}_0$ be the $\hat{Z}$-category of those 
Harish-Chandra modules having the generalized infinitesimal character of 
$\C$. 


\begin{ttt}
 There is a (uniquely determined up to isomorphism) 
 $\hat{Z}$-algebra $A$ and a finite set $X$ satisfying 
 \begin{itemize}
  \item ${\cal H}_0$ is equivalent, as $\hat{Z}$-category, to the 
   category $A$-{\bf fd} of finite dimensional $A$-modules, 
  \item $A$ is finitely generated over $\hat{Z}$, 
  \item $A/\rad(A)$ is isomorphic to the algebra of $\C$-valued functions 
   on $X$, 
  \item $A$ is selfopposite, 
  \item the global dimension of $A$ is $\dim G/K$, 
  \item the inclusion of $A$-{\bf fd} into $A$-{\bf mod} is compatible with 
   $\Ext$-calculus.
 \end{itemize}
\end{ttt}

{\bf Proof.} Let $n$ be a positive integer 
and ${\cal H}_n$ the full subcategory of ${\cal H}_0$ whose 
objects are killed by $\hat{I}^n$. This category comes with a selfduality 
and is equivalent to $A_n$-{\bf fd} where $A_n$ is a finite dimensional  
$(\hat{Z}/\hat{I}^n)$-algebra that is equipped with an anti-involution and 
$A_n/\rad(A_n)$ is isomorphic to the algebra of $\C$-valued 
functions on $X$. We have $A_{n+1}/\hat{I}^nA_{n+1}$ $=$ $A_n$ 
and $A_{n+1}/\hat{I}A_{n+1}$ $=$ $A_1$~; moreover the projections $A_{n+1}$ 
$\epi$ $A_n$ commute with the anti-involutions. Let $A$ be the limit of the 
$A_n$. Since $Z$ is a symmetric algebra over a finite dimensional vector 
space, it is noetherian, and so is 
$\hat{Z}$ by Proposition III.3.4.8 of Bourbaki in 
\cite{bac}, implying that $A$, being finitely generated over $\hat{Z}$, is 
itself noetherian. Section I.5.5 of Borel-Wallach \cite{bw} entails that 
the global dimension of $A$-{\bf fd} is $\dim G/K$ and that the 
inclusion of $A$-{\bf fd} into $A$-{\bf mod} is compatible with 
$\Ext$-calculus. 
The claim about the global dimension of $A$ now follows 
from statements 12 and 14 of Eilenberg in \cite{e1}. $\square$ \bb

\ind I hope there is always a partial ordering on $X$ which 
satisfies the assumptions of Theorem~\ref{Db}. 


\section{Example}

Let $\K$ be a commutative ring, let $z$ be an element of $\K$, 
let $A$ be the $\K$-algebra of the quiver 

$$
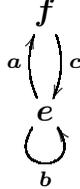
\xymatrix{
f\ar@/^/[d]^c\\
e\ar@/^/[u]^a\ar@(dl,dr)_b
}$$

modulo the relations

$$0=ab=bc=ac-zf=b^2+ca-ze.$$

\ind If $\K:=\C[[z]]$, where $z$ is an indeterminate, and 
${\cal H}_0$ is the category of those Harish-Chandra 
modules over $SL(2,\C)$ which have the generalized infinitesimal character 
of the trivial module, then ${\cal H}_0$ is equivalent to $A$-{\bf fd} 
(see Gelfand-Ponomarev \cite{gp}). \bb

\ind Put $X:=\{e,f\}$ with the ordering $e<f$. By Bergman's Diamond Lemma 
\cite{diamond} the set $\{e,f,a,b,c,b^2\}$ is a $\K$-basis of $A$, and 
we have 

$$eAe=\K\ e\oplus\K\ b\oplus\K\ b^2,\quad eAf=\K\ c,\quad fAe=\K\ a,\quad 
fAf=\K\ f,$$

$$M_e=Ae/AfAe=Ae/Aa=Ae/(\K\ a\oplus\K\ b^2),\ M_f=Af.$$

The sequence $M_f\mono Ae\epi M_e$, where the first arrow is right 
multiplication by $a$, being exact, $Ae$ belongs to 
$\langle M_x\rangle_{x\in X}$, and $A$ satisfies the assumptions of 
Theorem~\ref{Db}. \bb

\ind Denoting by $\dim R$ the global dimension of a ring $R$ and 
setting $B:=A/AfA$, we 
have $1+\dim\K\le\dim A\le2+\dim\K$. Indeed the isomorphism 
$B\simeq\K[t]$, where $t$ is an indeterminate, implies $\dim B=1+\dim\K$ 
(see for instance Theorem 4.3.7 in \cite{weibel}), 
and the claim follows from the proof of Lemma~\ref{extr} and the 
spectral sequence 

$$\Ext^p_B(V,\Ext^q_A(B,W))\then\Ext^{p+q}_A(V,W).$$ 

In particular the functor $i^*$ of Theorem~\ref{Db} exists if 
$\dim\K<\infty$.




\end{document}